\documentclass{article}
\usepackage{amsmath,amsfonts,amssymb,amsthm,amscd}
\title{Zero-groups and maximal tori}
\author{Alessandro Berarducci\footnote{URL:
www.dm.unipi.it/$\sim$berardu. Partially supported by Progetto
MIUR, Cofin 2004, Metodi di Logica in Algebra, Analisi e
Geometria. I thank Jean-Pierre Ressayre and the \'Equipe de
Logique of the University of Paris VII for their invitation and
kind hospitality in the period March 5 - April 5, 2005.}
\\Universit\`a di Pisa\\Largo Bruno Pontecorvo 5\\56127 Pisa}
\date{July 3, 2005 (revised Oct. 1)}
\DeclareMathOperator{\R}{\mathbb R}
\DeclareMathOperator{\N}{\mathbb N}

\DeclareMathOperator{\Z}{\mathbb Z}

\DeclareMathOperator{\M}{\bf M}

\DeclareMathOperator{\rank}{rank} \DeclareMathOperator{\im}{Im}
\DeclareMathOperator{\Card}{Card} \DeclareMathOperator{\Ker}{Ker}

\theoremstyle{plain}
\newtheorem{theorem}{Theorem}
\newtheorem{lemma}[theorem]{Lemma}
\newtheorem{proposition}[theorem]{Proposition}
\newtheorem{corollary}[theorem]{Corollary}

\newtheorem{claim}[theorem]{Claim}

\theoremstyle{definition}
\newtheorem{remark}[theorem]{Remark}
\newtheorem{definition}[theorem]{Definition}

\numberwithin{theorem}{section}
\begin{document}
\maketitle

\begin{abstract} We give a presentation of various results on
zero-groups in o-minimal structures together with some new
observations. In particular we prove that if $G$ is a definably
connected definably compact group in an o-minimal expansion of a real
closed field, then for any maximal definably connected abelian subgroup $T$ of
$G$, $G$ is the union of the conjugates of $T$. This can be seen as a
generalization of the classical theorem that a compact connected Lie
group is the union of the conjugates of any of its maximal
tori. \end{abstract}

\section{Introduction}
We consider groups definable in an o-minimal expansion $\M =
(M,<,+,\cdot, \ldots)$ of a real closed field $(M,<,+,\cdot)$.
Classical examples of such groups are the (real)algebraic subgroups of
the general linear group $GL(n,M)$. Identifying the algebraically
closed field $C = M[\sqrt{-1}]$ with $M^2$ in the standard way, we
also obtain all the algebraic subgroups of $GL(n,C)$. Less classical
examples of definable groups can be found in \cite{Strzebonski} or in
\cite{Peterzil-Steinhorn}. By \cite{Pillay1} each definable group can
be equipped with a unique group topology which makes it into a
definable manifold (see Definition \ref{definable-manifold}). A notion
of definable compactness for definable manifolds (and more generally
for definable spaces) can be introduced as in
\cite{Peterzil-Steinhorn}. (In the semialgebraic case the definably
compact spaces are the complete spaces of \cite{Knebusch}.) One has
also a notion of definable connectedness: a definable group is
definably connected if it has no proper definable subgroups of finite
index \cite{Pillay1}.  Unless $M=\R$, definable compactness (or
connectedness) does not imply compactness (or connectedness). In this
note we prove:

\begin{quote} {\bf Theorem \ref{conjugates2}}.
If $G$ is a definably connected definably compact group in an
o-minimal expansion $\M$ of a real closed field, then for any
maximal definable abelian subgroup $T$ of $G$, $G$ is the union of
the conjugates of $T$.\footnote{We learned from a referee that a
recent paper of M. Edmundo \cite{Edmundo-divisibility} contains another
proof of this result which uses the classification of definable
semisimple groups.}
\end{quote}

This can be seen as a generalization of the classical theorem that
a compact connected Lie group is the union of the conjugates of
any of its maximal tori. Indeed the classical theorem can be
deduced from our theorem taking $M = \R$, after some definability
considerations (Corollary \ref{tori}).

We show more generally that if $G$ is a definably connnected definably
compact group and $H<G$ is a definable subgroup of $G$ such that the
(o-minimal) Euler characteristic of $G/H$ is different from zero, then
$G$ is the union of the conjugates of $H$ (Theorem \ref{conjugates1}).

Under the hypothesis of Theorem \ref{conjugates2} it is not
difficult to prove that the ``Weyl group'' $W(G) := N_G(T)/T$ is
finite (Theorem \ref{weyl}). 

From Theorem \ref{conjugates2} and earlier work on definable abelian
groups one can obtain a result of M. Otero, also proved by M. Edmundo
in \cite{Edmundo-divisibility} by different methods, stating that
every definably compact definably connected group $G$ is divisible
(Corollary \ref{divisible}).

The proof of Theorem \ref{conjugates2} makes use of an o-minimal
version of Lefschetz fixed point theorem (see \cite{Edmundo2} or
\cite[Thm.  3.3]{Transfer}) and of the notion of $0$-group introduced
by Strzebonski in \cite{Strzebonski}. Strebonski generalized many
classical results on $p$-groups to the case of $0$-groups. It turns
out that, in the definably compact case, the $0$-groups are exactly
the definably connected abelian groups (Corollary
\ref{compact-zero-groups}).

Granted the appropriate background the proof of Theorem
\ref{conjugates2} is rather short, but we have taken this
opportunity to give an exposition, with bibliographical
references, of the relevant notions and results, and to make some
side observations.

\section{Euler characteristic}
We assume some familiarity with the basic notions of o-minimality
(see \cite{vdD}). Fix in the sequel an o-minimal structure $\M =
(M, <, +, \cdot, \ldots)$ expanding a real closed field (the dots
represent possible additional structure besides the field
structure). Although many of the results we use remain true also
for more general o-minimal structures, working over a field
simplifies some proofs and moreover the known proofs of the fixed
point theorem (Theorem \ref{fixed}) do make use of the field
structure. A subset $X\subset M^k$ is definable if it is (first
order) definable in $\M$ possibly with parameters. For instance,
if $\M$ is a real closed field $(M,<,+,\cdot)$ without additional
structure, the definable sets are exactly the semialgebraic sets.
We give to $M$ the topology generated by the open intervals and to
$M^n$ the product topology. To each definable set $X$ one can
attach two invariants: its dimension $\dim(X)\in \N$ and its
o-minimal Euler characteristic $E(X)\in \Z$.

\begin{definition} (see \cite{vdD}) The {\bf dimension} of $X$ is $\geq n$ iff
$X$ contains a subset definably homeomorphic to an open subset of
$M^n$. The (o-minimal) {\bf Euler characteristic} $E(X)$ is
defined as the number of even dimensional cells in $X$ minus the
number of the odd dimensional cells in $X$, relative to any given
cell decomposition. \end{definition}

The name ``Euler characteristic'' may be slightly misleading,
since $E(X)$ does not always coincide with the classical Euler
characteristic $\chi(X)$, even when the underlying o-minimal
structure $\M$ is based on the real numbers. For instance an open
interval $I \subset \R$ has $E(I) = -1$ and $\chi(I) = 1$. To
understand why this is the case one needs to remind that o-minimal
cells do not include the boundary, so an open interval is a odd
dimensional cell. This discrepancy is reflected in the fact that,
while $\chi$ is invariant under homotopies, $E$ satisfies instead
the following:

\begin{proposition} (see \cite[Ch. 4 (1.3),(2.4)]{vdD}) $\dim(X)$ and $E(X)$ are
invariant under definably bijections, not necessarily continuous.
\end{proposition}

Following \cite{Strzebonski} we can now define $\dim(G/H)$ and
$E(G/H)$ where $G$ is a definable group and $H<G$ is a definable
subgroup.

\begin{definition} A {\bf definable group} is a definable set $G\subseteq M^k$
together with a definable group operation. We do not require the
group operation to be continuous in the topology induced from
$M^k$. Let $G$ be a definable group, and let $H$ be a definable
subgroup, not necessarily normal. Let $G/H$ be the set of left
cosets of $H$. By ``definable choice'' (see \cite[Ch. 6 (1.2) p.
94]{vdD}) there is a definable function $\iota$ with domain $G$
such that $\iota (g) = \iota (h)$ iff $gH = hH$ (one can also
require $\iota (g) \in gH$). So whenever convenient we can
identify $G/H$ with the definable set $\iota (G)$ (identifying $gH
\in G/H$ with $\iota (g)$). Different choices of $\iota$ give rise
to definable bijections, so as a definable set $G/H$ is only
defined up to definable bijections. However, since the o-minimal
Euler characteristic and the dimension are invariant under
definable bijections, {\bf $E(G/H)$ and $\dim(G/H)$ are well
defined}. \end{definition}

It is possible to define, besides $E(X)$, another invariant $E'(X)$
which more closely resembles the classical Euler characteristic and is
invariant under definable homotopies. To this aim it suffices to
replace, in the classical definition $\chi(X) = \Sigma_i (-1)^i \rank
H_i(X)$, the classical homology group $H_i$ with the o-minimal
homology groups of Woerheide \cite{Woerheide} (which are naturally
isomorphic the classical ones when the o-minimal structure is based on
the reals, see \cite[Prop. 3.2]{Fundamental}). When $X$ is a closed
and bounded subset of $M^k$ one can easily prove using the
triangulation theorem that $E(X) = E'(X)$ (see
\cite[p. 788]{Transfer}), so in particular if $M=\R$ and $X\subset
\R^k$ is compact, then $E(X) = E'(X) = \chi(X)$. Thus $E(X)$ is a tool
which allows one to use combinatorial arguments typical of finite
groups theory (thanks to its invariance under definable bijections),
while at the same time permitting to draw conclusions of topological
nature (thanks to the fact that it coincides with $E'(X)$ when $X$ is
closed and bounded). We will implicitly use these facts in the
sequel. In fact the results relying on the fixed point theorem depend
on the use of $E'$.

\section{Definable spaces}
The notion of definable space is discussed in \cite{vdD}. We give
here an apparently weaker definition that is equivalent up to
isomorphisms.

\begin{definition} \label{definable-manifold} Let $X$ be a definable set and let
$\tau$ be a topology on $X$. We say that $(X, \tau)$ is a {\bf
definable space} if there is a finite cover of $X$ by definable
sets $U_1, \ldots, U_k$, open in $\tau$, and a natural number $n$,
such that for each $U_i$ there is a definable function $\varphi_i
\colon U_i \to V_i$ which is a homeomorphism between $U_i$ (with
the $\tau$ topology) and a subset $V_i$ of $M^n$ (where $M$ has
the topology generated by the open intervals, and $M^n$ has the
product topology). If moreover the $V_i$ are open in $M^n$, we say
that $(X, \tau)$ is a {\bf definable manifold} of dimension $n$.
The collection of the ``{\bf local charts}'' $(U_i, \varphi_i, n)$
is called an {\bf atlas} of $(X,\tau)$. So each definable space or
manifold has a finite atlas.

The easiest example of definable space is a {\bf definable
subspace of $M^n$}, namely a definable set $X\subset M^n$ with the
topology inherited by the ambient space $M^n$. A definable space
is {\bf affine} if it can be definably embedded in $M^n$ for some
$n$ (namely if it is definably homeomorphic to a definable subset
of $M^n$ considered as a definable subspace of $M^n$). \end{definition}

We will later need:

\begin{theorem} (\cite{Robson}, see also \cite[Ch. 10 (1.8)]{vdD})
\label{Robson} A necessary and sufficient condition for a
definable space to be affine is that it is {\bf regular}, namely
its points are closed and for every point $p$ and closed set $C$
there are open neighbourhoods of $p$ and $C$ which are disjoint.
\end{theorem}

Since each open set is a union of definable open sets, it is easy to
see that a definable space is regular if and only if it is {\bf
definably regular}, namely its points are closed and for every point
$p$ and definable closed set $C$ there are definable open
neighbourhoods of $p$ and $C$ which are disjoint.

\begin{lemma} \label{submanifold} Let $(X, \tau)$ be a definable manifold of
dimension $n$ and let $Y\subset X$ be a subset of $X$ of dimension
$m$, equipped with the subspace topology.  Then $Y$ is a definable
(sub)manifold of $X$ if and only if every point $p\in Y$ has a
definable open neighbourhood $O\subset Y$ definably homeomorphic
to an open subset of $M^m$. \end{lemma}

So although definable manifolds are required to have a finite
atlas, for submanifolds it is not necessary to require this
finiteness condition as it is always automatically ensured.

\begin{proof} The lemma was proved in \cite[Prop. 4.2]{Fundamental} in the
case when $(X,\tau) =M^n$. The general case can be reduced to this
special case using the fact that $X$ is covered by finitely many
open sets definably homeomorphic to open subsets of $M^n$. \end{proof}

Our main source of definable spaces, besides the subspaces of $M^n$,
are the definable groups. Recall that for a definable group
$G\subseteq M^k$ we do not require the group operation to be
continuous in the topology induced from the ambient space $M^k$.  The
good topology on definable groups is not the topology of the ambient
space but the one given by the following:

\begin{theorem} (\cite[Prop. 2.5]{Pillay1}) \label{group-topology} If $G$ is a
definable group, then there is a (unique) topology $\tau$ on $G$,
called the {\bf definable manifold topology} of $G$, such that:
\begin{enumerate}
\item $(G, \tau)$ is a topological group (i.e. multiplication and
inversion are continuous);
\item $(G, \tau)$ is a definable
manifold of dimension $n$, where $n$ is the (o-minimal) dimension
of $G$.
\end{enumerate}
\end{theorem}

The uniqueness of the definable group topology follows from:

\begin{lemma} (\cite[Lemma 1.11]{PPS1}) \label{uniqueness} If $f \colon H
\to G$ is a definable group homomorphism between definable groups
equipped with definable manifold topologies, then $f$ is
continuous. \end{lemma}

For definable spaces we have the following notion of definable
compactness:

\begin{definition} (\cite{Peterzil-Steinhorn}) An Hausdorff definable space
$(X,\tau)$ is {\bf definably compact} iff for every definable
function $f\colon I \to X$, where $I = (a,b) \subset M$ is an open
interval, $\lim_{t \to b^-} f(t)$ exists in $(X,\tau)$. \end{definition}

Since by o-minimality every definable function $f\colon I \to
(X,\tau)$ is piecewise continuous, without loss of generality we
can assume $f$ continuous.

\begin{proposition} (\cite[Thm. 2.1]{Peterzil-Steinhorn}) Let $X\subset M^n$ be a
definable set with the topology induced by the ambient space
$M^n$. Then $X$ is definably compact iff it is closed and
bounded.\end{proposition}

The closed and bounded subsets of $M^n$ need not be compact if
$M\neq \R$, but they behave in many respects like compact sets
within the definable category. For instance the image of a closed
and bounded set under a continuous definable function is closed
and bounded (see \cite[Ch. 6 (1.9) p. 95]{vdD} for o-minimal
structures expanding an ordered abelian group, and
\cite{Peterzil-Steinhorn} for arbitrary o-minimal structures).

\begin{lemma} \label{closed} Let $f\colon (X, \tau) \to (Y, \mu)$ be a
definable continuous surjective function between definable spaces
and assume that $(X, \tau)$ is definably compact. Then $(Y, \mu)$
is definably compact and $f$ is a closed map. \end{lemma}

\begin{proof} Let $\gamma \colon I \to Y$ be a definable function, $I =
(a,b)$. By definable choice $\gamma$ can be lifted to a definable
function $\sigma \colon I \to X$ with $f \circ \sigma = \gamma$.
Since $\sigma$ has a limit in $X$ (and $\gamma$ and $\sigma$ are
piecewise continuous by o-minimality), $f \circ \sigma$ has a
limit in $Y$. So $Y$ is definably compact. Since a closed subset
of a definably compact set is definably compact, the same argument
shows that $f$ is a closed map. \end{proof}

\begin{theorem} (\cite[Cor. 2.8]{Pillay1}) Let $G$ be a definable group and
let $H<G$ be a definable subgroup. Then $H$ is closed in the
definable group topology of $G$.  \end{theorem}

Moreover we have:

\begin{theorem} (\cite[Lemma 2.6]{Peterzil-Starchenko}) Let $G$ be a definable group and $H<G$ be a definable
subgroup. Then the definable group topology on $H$ coincides with
the topology as a subspace of $G$ \end{theorem}

The proof in \cite{Peterzil-Starchenko} is in term of generic
elements. The next lemma yields a proof based on Lemma
\ref{submanifold}.

\begin{lemma} \label{fH} Let $f\colon H \to G$ be an injective definable
morphism of definable groups. Then $f$ is a homeomorphism from $H$
to $f(H)$, where $H$ and $G$ have the definable group topology and
$f(H)\subseteq G$ has the subspace topology. \end{lemma}

\begin{proof} $f \colon H\to G$ is continuous by Lemma \ref{uniqueness}. Let
$K = f(H) \subseteq G$ and let $m = \dim H = \dim K$.

By the cell decomposition theorem $K$ contains an open subset $V
\subset K$ (in the topology inherited from $G$) definably
homeomorphic to an open subset of $M^m$ (take a cell of dimension
$m$ in the intersection of $K$ with a local chart of $G$).

Now $K$ is a subgroup of $G$ and since it has the subspace
topology it is also a topological subgroup. This implies that the
translates of $V$ in $K$ are definably homeomorphic to each other,
and since they cover $K$, we get that every point of $K$ has a
neighbourhood definably homemomorphic to an open subset of $M^m$.
The number of such translates need not be finite, but nevertheless
by Lemma \ref{submanifold} we obtain that $K$ is a definable
submanifold of $G$. We can then apply Lemma \ref{uniqueness} to
conclude that $f^{-1}\colon K \to H$ is continuous, thus finishing
the proof. \end{proof}

If $H$ is definably compact Lemma \ref{fH} has a shorter proof.
Indeed in this case the continuous map $f \colon H\to G$ is a
closed map (by Lemma \ref{closed}), hence an homeomorphism onto
its image.

We will later need:

\begin{definition} \label{connected} A definable space $(X,\tau)$ is {\bf
definably connected} if and only if it has no definable proper
non-empty clopen subset. \end{definition}

\begin{remark} Each definable space $(X,\tau)$ is a finite union of maximal
definably connected subsets, called its definably connected
components.\end{remark}

\begin{proof} Immediate from the cell decomposition theorem when $X$ is a
subspace of $M^k$. The general case follows working in the local
charts. \end{proof}

\begin{remark} \label{path-connected} A definable space $(X,\tau)$ is
definably connected if and only if it is definably path connected,
i.e. any two points can be joined by a definable continuous path.
(See \cite[Ch. 6, Prop. 3.2, p. 100]{vdD} for the case of
subspaces of $M^k$.) \end{remark}

\section{Homogeneous spaces}
Besides definable groups, another source of definable manifolds
are the definable homogeneous spaces, namely the definable sets on
which a definable group acts transitively by a definable action.

\begin{theorem} (\cite[Thm. 2.11]{PPS1}) \label{action1} Let $G$ be a
definable group, and let $\alpha \colon G \times V \to V$ be a
definable transitive action on the definable set $V$.\footnote{We
recall that the action $\alpha$ is transitive if for every
$v,v'\in V$ there is $g\in G$ with $\alpha(g, v) = v'$.} There is
a topology $\tau$ on $V$ such that:
\begin{enumerate}
\item
 $(V, \tau)$ is a definable manifold;
\item the
action $\alpha \colon G \times V \to V$  is continuous (where $G$
has the definable manifold topology).
\end{enumerate}
\end{theorem}

\begin{corollary} (\cite[Cor. 2.14]{PPS1}) \label{action} Let $G$ be a definable
group, and let $H<G$ be a definable subgroup, not necessarily
normal. There is a topology $\tau$ on $G/H$ such that:
\begin{enumerate}
\item
 $(G/H, \tau)$ is a definable manifold;
\item the natural
action $L\colon G \times G/H \to G/H$ given by left multiplication
$L(g_1, g_2H) = g_1g_2H$ is continuous (where $G$ has the
definable manifold topology).
\end{enumerate}
\end{corollary}

\begin{theorem} \label{coincide} Let $G$ be a definable group and let $H<G$ be
a definable subgroup. The following topologies on $G/H$ coincide:
\begin{enumerate}
\item The topology $\tau$ on $G/H$ given by Theorem \ref{action};
\item The quotient topology $\nu$ on $G/H$ (where $G$ has the definable manifold topology);
\item If $H$ is normal in $G$, the group $G/H$ has its own definable manifold topology given by Theorem
\ref{group-topology}, which also coincides with the quotient
topology.
\end{enumerate}
\end{theorem}

Thanks to this result we can speak of the {\bf definable manifold
topology} on $G/H$ without ambiguity.

\begin{proof} The projection $\pi\colon G \to (G/H, \tau)$ is continuous
because $\pi(g) = gH = L(g,1H)$ and $L$ is continuous. Since the
kernel is $H$, passing to the quotient we have an induced
continuous map $id\colon (G/H, \nu) \to (G/H, \tau)$ where $\nu$
is the quotient topology and $id$ is the identity. To prove that
$\nu = \tau$ it remains to show that $id\colon (G/H, \nu) \to
(G/H, \tau)$ is an open map (hence an homeomorphism). By
definition of the quotient topology this is equivalent to say that
$\pi \colon G \to (G/H, \tau)$ is an open map, and this is ensured
by Theorem \ref{openmap} below. Finally if $H\lhd G$, we know that
$G/H$ has a unique group topology which is at the same time a
definable manifold topology. Since the quotient topology $\nu$ on
$G/H$ is a group topology, and we have just shown that it
coincides with the definable manifold topology $\tau$, the proof
is finished.  \end{proof}

It remains to prove:

\begin{theorem} \label{openmap} Let $\alpha \colon G \times V\to V$ be a
definable transitive action and let $\tau$ be the definable
manifold topology on $V$ given by Theorem \ref{action1}. Let $v\in
V$. The continuous map $\gamma \colon G \to (V,\tau)$ defined by
$\gamma (g) = \alpha(g, v)$, is an open map. \end{theorem}

\begin{proof} We say that $\gamma$ is open at $p\in G$ if for every open
neighborhood $O$ of $p$, $\gamma(O)$ contains $\gamma(p)$ in its
interior. Clearly if a $\gamma$ is open at every point then it is
an open map. Using the transitivity and the continuity of the
action of $G$ (on both $G$ and $V$), it is easy to see that if
$\gamma$ is open at some point then it is open at every point. So
it suffices to prove that $\gamma$ is open at some point. We will
work on local charts and we will use the following easy
consequence of the trivialization theorem (see \cite[Ch. 9 (1.2)
p. 142]{vdD} for the statement of the trivialization theorem):

\begin{claim} Let $X,Y \subseteq M^k$ be definable sets with the topology
induced by $M^k$ and let $f \colon X \to Y$ be a continuous onto
map. Then $f$ is open at some point of $X$.
\end{claim}

To prove the claim recall that $f\colon X \to Y$ is definably trivial
if there is a definable homeomorphism $\sigma \colon X \to Y \times F$
such that $f = \pi_1 \circ \sigma$, where $\pi_1 \colon Y\times F \to
Y$ is the projection. Clearly if $f$ is definably trivial it is an
open map, since so are the projections.  In the general case, by the
trivialization theorem and the surjectivity of $f$ we can find an open
set $O$ of $Y$ such that the restriction of $f$ to $f^{-1}(O)$ is a
definably trivial map onto $O$, and the claim follows.

To finish the proof consider a cell decomposition of $G$
compatible with the open sets of a finite atlas of $G$ and with
the $\gamma$-preimages of the open sets of a finite atlas of $V$.
At least one cell of this decomposition of $G$ is mapped by
$\gamma$ into a subset of $V$ containing a non-empty open set
$P\subset V$. Let $X=\gamma^{-1}(P) \subseteq G$. The restriction
of $\gamma$ to $X$ is a surjective continuous map $\gamma_{|X}
\colon X \to P$. Moreover by our costruction $X$ is contained in a
single chart of $G$, and $P$ is contained in a single chart of
$V$. So using the claim we can conclude that $\gamma_{|X}$ is open
at some point $p$ of $X$, so a fortiori $\gamma$ is open at $p$.
\end{proof}

\begin{remark} If $G$ is definably compact the proof of Theorem
\ref{coincide} can be simplified. In fact by Lemma \ref{closed}
and the continuity of $\pi$, if $G$ is definably compact, also
$(G/H, \tau)$ is definably compact and the continuous injective
map $id\colon (G/H, \nu) \to (G/H, \tau)$ is a closed map, hence
an homeomorphism. \end{remark}

\section{Zero-groups}

\begin{definition} $G$ is a {\bf $0$-group} if it is a definable group and for
every proper definable subgroup $H<G$ we have $E(G/H) = 0$ (we
allow $H$ to be the trivial group). \end{definition}

Note that the trivial group is a $0$-group since it has no proper
subgroups.

\begin{theorem} (\cite[Prop. 2.12]{Pillay1}) Let $G$ be a definable group. Put
on $G$ the definable manifold topology. Let $G^0\subset G$ be the
definably connected component of the identity of $G$. Then $G^0$
is a normal subgroup of $G$ and it is the smallest definable
subgroup of $G$ of finite index. Thus $G$ is definably connected
if and only if it has no proper definable subgroups of finite
index. \end{theorem}

The definition of definable connectedness given above (Definition
\ref{connected}) is thus consistent with the following:

\begin{definition} A definable group $G$ is {\bf definably connected} iff $G$ has
no definable proper subgroup of finite index. \end{definition}

\begin{remark} If $G$ is a $0$-group, then $G$ is definably connected. \end{remark}

\begin{proof} Let $G^0 <G$ be the definably connected component of the
identity. Then $[G:G^0]$ is finite and therefore $E(G/G^0) =
[G:G^0] \neq 0$. If $G^0 \neq G$ this would contradict the
definition of $0$-group. \end{proof}

\begin{theorem} (\cite[Cor. 5.7]{Strzebonski}) If $G$ is a $0$-groups, then
$G$ is abelian. \end{theorem}

\begin{theorem} (\cite{Edmundo2} or \cite[Thm. 3.3]{Transfer}) \label{Ezero}
If $G$ is a definably compact infinite group, $E(G) = 0$. \end{theorem}

\begin{corollary} \label{compact-zero-groups} Let $G$ be a definably compact
group. Then $G$ is a $0$-group if and only if it is abelian and
definably connected. \end{corollary}

\begin{proof} Suppose $G$ is definably compact, abelian and definably
connected. We must prove that $G$ is a $0$-group. Let $H<G$ be a
proper subgroup (if $H$ does not exist, $G$ is the trivial group,
which is a $0$-group). Since $G$ is abelian $H\lhd G$. The group
$G/H$ is definably compact since so is $G$, and it is infinite
since $G$ is definably connected. So $E(G/H)=0$ by Theorem
\ref{Ezero}. \end{proof}

Strzebonki gave an example of a definable (even real
semialgebraic) $0$-group which is not definably compact \cite[Ex.
5.3]{Strzebonski}. The $0$-groups which are definably compact are
exactly the Strzebonski tori defined below:

\begin{definition} $G$ is a Strzebonski torus if and only if $G$ and all
definably connected subgroups of $G$ are $0$-groups. \end{definition}

Since $0$-groups are definably connected, this is equivalent to
the definition given by Strzebonski in \cite{Strzebonski}: $G$ is
a (Strzebonski) torus if $G$ is a $0$-group and every $H<G$
contains a $0$-group $K<H$ with $[H:K]$ finite.

\begin{proposition} Let $G$ be a definable group. Then $G$ is a Strzebonski torus
iff $G$ is definably compact abelian and definably connected. \end{proposition}

\begin{proof} Assume $G$ is a Strzebonski torus. Then in particular is a
$0$-group, so it is abelian and definably connected. It remains to
show that $G$ is definably compact. We reason by contradiction
using a result of Peterzil and Steinhorn
\cite{Peterzil-Steinhorn}: if a definably group $G$ is not
definably compact, then it contains a definable one-dimensional
torsion free subgroup $H<G$. We can assume that $H$ is definably
connected (as otherwise replace it with its component at the
identity). Since groups with $E = 0$ have elements of every prime
order (by \cite{Strzebonski}), it follows that $E(H)\neq 0$, so
$G$ is not a Strzebonski torus.

For the opposite direction suppose $G$ is definably compact abelian and definably
connected. We have already proved that $G$ is a $0$-group. To show
that it is a torus we must show that if $H<G$ is definably
connected, then $H$ is a $0$-group. This is clear since $H$
satisfies the same assumptions used to show that $G$ is a
$0$-group. \end{proof}

\section{$0$-Sylow subgroups and maximal tori}

\begin{theorem} \label{fixed3} Let $G$ be a definably compact group and let
$H<G$ be a definable subgroup. If $E(G/H) \neq 0$ and $f \colon
G/H \to G/H$ is a definable continuous map definably homotopic to
the identity (with respect to the definable manifold topology of
$G/H$), then $f$ has a fixed point. \end{theorem}

\begin{proof} Let us first observe:

\begin{claim} $G/H$ is a regular topological space
(we do not need definable compactness here).
\end{claim}

Indeed each definable manifold, although it does not need to be
Hausdorff, it is certainly $T_1$ (i.e. its points are closed). So
$G$ with the definable manifold topology is $T_1$. Now we use the
fact that if $G$ is a $T_1$ topological group and $H<G$ is a
closed subgroup, then $G/H$ is regular (see \cite[Ch. 1, Thm.
1.6]{Mimura}).

By Theorem \ref{Robson} we can thus identify $G/H$ with a
definable submanifold of $M^k$. For such manifolds a definable
version of the singular homology groups has been developed by
Woerheide (\cite{Woerheide}), and a notion of orientability can
then be defined as in \cite{Transfer} in terms of the definable
homology groups. A routine argument already used in
\cite{Transfer} (for $G$ instead of $G/H$) shows:

\begin{claim} $G/H$ is definably orientable. \end{claim}

The idea is that (as for classical homogeneous spaces) an
orientation on $G/H$ is obtained by choosing a local orientation
at a point and extending it consistently to the whole space $G/H$
using the transitivity of the action of $G$ on $G/H$.

To conclude the proof of Theorem \ref{fixed3} it now suffices to
invoke Theorem \ref{fixed} below. \end{proof}

\begin{theorem} (\cite{Edmundo2} or \cite[Thm. 3.3]{Transfer}) \label{fixed}
Let $X\subseteq M^k$ be a definable set. Suppose that $X$, with
the subspace topology from $M^k$, is a definably compact definably
oriented definable manifold with $E(X) \neq 0$. Let $f \colon X
\to X$ be a definable continuous map definably homotopic to the
identity. Then $f$ has a fixed point. \end{theorem}

We can now prove:

\begin{theorem} \label{conjugates1} Let $G$ be a definably connected definably
compact group and let $H<G$ be a definable subgroup with $E(G/H)
\neq 0$. Then $G = \bigcup_{g \in G} gHg^{-1}$. \end{theorem}

\begin{proof} Let $g\in G$. We want to prove that $g \in \bigcup_{x\in G}
xHx^{-1}$. Consider the map $L_g \colon G/H \to G/H$ defined by
$L_g(xH) = \alpha(g,xH) = gxH$. We claim that $L_g$ is definably
homotopic to the identity. Granted the claim, by Theorem
\ref{fixed3} $L_g$ has a fixed point $xH\in G/H$. So $gxH = xH$
and therefore $g\in xHx^{-1}$. To prove the claim recall, by
\ref{path-connected}, that $G$ is definably path connected. So
there is a definable continuous function $t \mapsto g_t$, $t\in
[0,1]^M$, with $g_0 = 1, g_1 = g$. Hence $t \mapsto L_{g_t}$ is a
definable homotopy between $id = L_{g_0}$ and $L_g = L_{g_1}$. \end{proof}

To apply the above results to $0$-subgroups we recall some results
of Strzebonski.

\begin{definition} Let $G$ be a definable group. A $0$-group $H<G$ is called
$0$-Sylow if it is a maximal $0$-subgroup of $G$. \end{definition}

\begin{remark} \label{bijections} (\cite[Thm. 2.14]{Strzebonski}) if $K<H<G$
are definable groups (no normality assumptions), there is a
definable bijection from $G/H \times H/K$ to $G/K$. So
$E(G/H)E(H/K) = E(G/K)$ and $\dim(G/H) + \dim(H/K) = \dim(G/K)$.
\end{remark}

\begin{proposition} (see \cite[Rem. 2.18]{Strzebonski}) If $H,G$ are $0$-groups
and $H$ is a proper subgroup of $G$, then $\dim(H) < \dim(G)$.
\end{proposition}

\begin{proof} By Remark \ref{bijections} $\dim(G) = \dim(G/H)+\dim(H)$.
Since $G$ is a $0$-group, $E(G/H) = 0$. So $G/H$ is infinite
(because for a finite set $X$, $E(X) = \Card(X)$). But then
$\dim(G/H)>0$ (since infinite definable sets have positive
dimension) and the result follows.\end{proof}

\begin{corollary} (\cite[Rem. 2.18]{Strzebonski}) Every $0$-subgroup is
contained in a $0$-Sylow. \end{corollary}

\begin{theorem} (\cite{Strzebonski}) \label{0Sylow} Let $G$ be a definable
group and let $H<G$ be a definable subgroup which is a $0$-group.
\begin{enumerate}
\item If $E(G/H)=0$,
  then there is a $0$-subgroup $K$ of $G$ with $H<K<N_G(H)$ and $K\neq
  H$.
\item If $E(G/H) \neq 0$, then $H$ is $0$-Sylow.
\item Thus $E(G/H) \neq 0$ if and only if $H$ is
$0$-Sylow.
\end{enumerate}
\end{theorem} \begin{proof}  Part 1. is (\cite[Thm. 2.14]{Strzebonski}). Part 2.
follows immediately from the definitions and Remark
\ref{bijections}. \end{proof}

\begin{theorem} \label{conjugates} Let $G$ be definably compact and definably
connected. Let $T<G$ be a $0$-Sylow of $G$. Then $G =
\bigcup_{x\in G} xTx^{-1}$. (Moreover by \cite[Thm.
2.21]{Strzebonski} each two $0$-Sylow subgroups are conjugates.)
\end{theorem}

\begin{proof} By Theorem \ref{conjugates1}. \end{proof}

Since in the definably compact case the $0$-groups are the
definably connected abelian groups the above theorem has the
following equivalent formulation:

\begin{theorem} \label{conjugates2} If $G$ is a definably connected definably
compact group, then for any maximal definably connected abelian subgroup $T$
of $G$, $G$ is the union of the conjugates of $T$. \end{theorem}

If the o-minimal structure is an expansion of $\R$ we obtain the
following classical result:

\begin{corollary} \label{tori} If $G$ is a compact connected Lie group, then for
any maximal abelian connected closed subgroup of $H<G$, $G$ is the union
of the conjugates of $H$. \end{corollary}

\begin{proof} The only thing to observe is that we do not need any
definability assumptions. The definability comes for free since
any compact Lie group $G$ is isomorphic to a compact subgroup $K$
of $GL(n,\R)$ for some $n$ (see \cite[Ch. 3, Thm. 4.1, p.
136]{Brocker}) and any such $K$ is a (real)algebraic subgroup of
$GL(n,\R)$ \cite[Prop. 2, p. 230] {Chevalley}, hence it is definable
in the o-minimal structure $(\R,<,+,\cdot)$. \end{proof}

Let $N_G(T)$ be the normalizer of $T$ in $G$. In the hypothesis of
Theorem \ref{conjugates2} we have:

\begin{theorem} \label{weyl} The ``Weyl group'' $W(G) := N_G(T)/T$ is finite.
\end{theorem}

\begin{proof} Suppose $W = W(G)$ is infinite. Then so is its definably
connected component $W^0<W$. Since $G$ is definably compact, so is
$W^0$. Hence $E(W^0)=0$. Now $E(W)=E(W/W^0)E(W^0) = 0$. By
definition this means $E(N_G(T)/T)=0$. Hence, by Theorem \ref{0Sylow}, $T$ is
not a zero-Sylow, a contradiction.
\end{proof}

A compact connected abelian Lie group is isomorphic to a torus,
namely a product of finitely many copies of $\R/\Z$. This suggests
that a definably compact definably connected abelian definable
group $G$ of (o-minimal) dimension $n$ should resemble an
$n$-dimensional torus ${(\R/\Z)}^n$. This analogy is partially
confirmed in \cite{Edmundo-Otero}, where it is shown that such a
group $G$ has the same torsion of an $n$-dimensional torus.
Further information on $G$ comes by considering quotients.
Assuming the o-minimal structure sufficiently saturated, in
\cite{BOPP} it is shown that $G$ has a smallest type-definable
subgroup $G^{00} < G$ of bounded index (which is moreover normal
and divisible). By \cite[Theorem 8.1]{HPP} $G^{00}$ is torsion
free and $G/G^{00}$ isomorphic to the torus $(\R/\Z)^n$.

However the analogy with tori has its limitations: for instance
there are definably compact definably connected abelian definable
group $G$ of dimension $n>1$ without any proper infinite definable
subgroup (see \cite{Peterzil-Steinhorn}). Hence, given a definably
connected group $G$ of dimension $>1$, we have no control on how
big a minimal $0$-subgroup $H$ is: it is not necessarily the case
that $\dim(H) = 1$ ($H$ could be $G$ itself).

\begin{corollary} (\cite{Otero-divisibility}, \cite{Edmundo-divisibility})
\label{divisible} If $G$
is a definably connected definably compact group, then $G$ is
divisible.  \end{corollary}

\begin{proof} By Corollary \ref{conjugates2} we can reduce to the case where
$G$ is abelian. The divisibility of definably connected abelian
groups follows from the results of \cite{Strzebonski}, as observed
in \cite{Edmundo-Otero}. In fact for $k \in \N$ define $p_k \colon
H \to H$ by $p_k(x) = x^k$. Since $H$ is abelian $p_k$ is a
homomorphism. By \cite{Strzebonski} the torsion subgroups $\Ker
p_k$ is finite, and hence $\dim H = \dim \im p_k$. Since $H$ is
definably connected $H = \im p_k$, hence $H$ is divisible. \end{proof}

\end{document}